\documentclass[12pt]{amsart}
\usepackage{graphicx}
\usepackage{caption}
\usepackage{subcaption}
\usepackage{color}
\usepackage{amsmath,amsfonts}
\usepackage{algorithm}
\usepackage{algorithmic}
\usepackage{booktabs}
\newtheorem{thm}{Theorem}[section]

\theoremstyle{definition}

\theoremstyle{remark}

\newtheorem{remark}[thm]{Remark}

\begin{document}

%%
%% The title of the paper goes here.  Edit to your title.
%%

%\title[Robust FSI Solver]{A note on the robust preconditioner for the monolithic fluid-structure interaction system}
\title[Robust FSI Solver]{A note on robust preconditioners 
for monolithic fluid-structure interaction systems 
of finite element equations}
%%
%% Now edit the following to give your name and address:
%% 
\author{Ulrich Langer}
\address{Johann Radon Institute for Computational and Applied Mathematics (RICAM), Austrian Academy of Sciences, Altenberger Strasse 69, A-4040 Linz, Austria}
\email{ulrich.langer@ricam.oeaw.ac.at}
\urladdr{http://www.ricam.oeaw.ac.at/people/u.langer/}% Delete if not wanted.

%%
%% If there is another author uncomment and edit the following.
%%

\author{Huidong Yang}
\address{Johann Radon Institute for Computational and Applied Mathematics (RICAM), Austrian Academy of Sciences, Altenberger Strasse 69, A-4040 Linz, Austria}
\email{huidong.yang@oeaw.ac.at}
\urladdr{http://people.ricam.oeaw.ac.at/h.yang/}

%%
%% If there are three of more authors they are added in the obvious
%% way. 
%%

%%%
%%% The following is for the abstract.  The abstract is optional and
%%% if not used just delete, or comment out, the following.
%%%

\begin{abstract}
  In this note, we consider preconditioned Krylov 
  subspace methods for discrete fluid-structure interaction 
  problems with a nonlinear hyperelastic 
  %model, 
  material model and covering a large range of flows, e.g, water, 
  blood, and air with highly varying density. 
  Based on the complete $LDU$ factorization of 
  the coupled system matrix, 
  the preconditioner is constructed 
  in form of 
  %as 
  $\hat{L}\hat{D}\hat{U}$, where $\hat{L}$, $\hat{D}$ and $\hat{U}$ 
  are 
  %the 
  proper 
  %approximation 
  approximations 
  to $L$, $D$ and $U$, respectively. The inverse of 
  %the 
  the corresponding 
  Schur complement is approximated by applying one 
  %iteration 
  cycle of a special class of algebraic multigrid methods 
  to the perturbed fluid sub-problem, that is obtained by modifying 
  corresponding entries in the original fluid matrix with an explicitly 
  constructed approximation of the exact 
  perturbation coming from the sparse matrix-matrix multiplications. 
\end{abstract}

%%
%%  LaTeX will not make the title for the paper unless told to do so.
%%  This is done by uncommenting the following.
%%

\maketitle

%%
%% LaTeX can automatically make a table of contents.  This is done by
%% uncommenting the following:
%%

%\tableofcontents

%%
%% A Theorem is stated by
%%

%\begin{thm} The square of any real number is non-negative.
%\end{thm}

%%
%% Its proof is set off by
%% 

%\begin{proof}
%\end{proof}

%%
%% A new section is started as follows:
%%

%%%%%%%%%%%%%%%%%%%%%%%%%%%%%%%%%%%%%%%%%%%%%%%%%%%%%%%%%%%%%%%%%%%%%%
\section{Motivation}\label{sec:intro}
During the past years, robust and efficient monolithic 
fluid-structure interaction (FSI) solvers attract 
%lost 
a lot 
of interests from many researchers; see, 
e.g., \cite{NME:NME3001, Razzaq20121156, STurek06, Heil20041, 
Badia20084216,EPFL11,ABXCC10}, 
that are mainly based on 
%the 
algebraic multigrid 
%(AMG, see \cite{JWRKS87}) 
(AMG) \cite{JWRKS87} 
%or 
, 
geometry multigrid 
%(GMG, see \cite{WH03}), 
(GMG) \cite{WH03}, 
preconditioned Krylov subspace 
%(see \cite{YS03:00}) 
\cite{YS03:00} 
and domain decomposition 
%(DD, see \cite{QuarteroniVali:1999,ATOW05}) 
(DD) \cite{QuarteroniVali:1999,ATOW05} 
methods. In our 
previous work \cite{UY14}, we implemented 
FSI monolithic AMG 
%a monolithic AMG 
%and a variant of the W-cycle 
%(a recursive Krylov-based multigrid cycle, 
%somehow related to the algebraic multilevel method, 
%see, e.g., \cite{AO89I,AO90,UJ:91, KJMS13, PSV08:00, NLA:NLA542, PSV14}) 
solvers with the W-cycle and with a variant of the W-cycle, 
i.e., a recursive Krylov-based multigrid cycle, 
somehow related to the algebraic multilevel method, 
see, e.g., \cite{AO89I,AO90,UJ:91, KJMS13, PSV08:00, NLA:NLA542, PSV14} 
%.FSI solvers 
, as well as their corresponding preconditioners 
for the coupled FSI problem using the AMG preconditioners 
\cite{FK98:00, UH:02, WM04:00, CVULHY1300} 
for each sub-problem in the smoothing steps. 
In addition, we also considered 
%therein 
the preconditioned GMRES method (see \cite{Saad86}), using a class of 
block-wise Gauss-Seidel 
%as 
type 
preconditioners (see the earlier work in \cite{NME:NME3001}), 
that are 
based on the aforementioned AMG methods for the sub-problems. 
As well known, such block-wise Gauss-Seidel preconditioned 
Krylov subspace methods may lose the robustness with 
respect to the mesh size, i.e., the iteration numbers 
for solving the coupled FSI system nearly double, when the mesh 
size halves; see the numerical results in our previous work \cite{UY14}. 
%Herein
In this work, 
we further obverse, that the methods are not 
robust with respect to the varying fluid density. For an illustration, we have 
tested the methods on the numerical example given in Section \ref{sec:num}, 
where the fluid density is varying (from water to air flow): 
$\rho_f\in\{1.1, 0.11, 0.011, 0.0011\}$ g/cm$^3$. 
The preconditioners employed in the 
preconditioned GMRES method are: The block diagonal ($\tilde{P}_D$), the 
block lower triangular ($\tilde{P}_L$), the block upper triangular ($\tilde{P}_L$), 
the SSOR ($\tilde{P}_{SSOR}$) and the $ILU(0)$ ($\tilde{P}_{ILU}$); see 
the definition in \cite{UY14}. The 
number of iterations (\#it) for solving the linearized FSI system on the coarse mesh using 
the time step size $\Delta t = 0.125$ ms is displayed in Table \ref{tab:unrob}.
\begin{table}[ht!]
\centering
\begin{tabular}{cccccc}
\toprule
&\multicolumn{5}{c}{\#it}  \\
 \cmidrule(r){2-6} 
 $\rho_f$
 &{1.1} &{0.11}
&{0.011}&{0.0011}
 &{0.00011} \\ 
\midrule
$\tilde{P}_D$                &$57$ &$185$ &$>250$   &$>250$ &$>250$ \\
$\tilde{P}_L$                &$38$ &$68$ &$89$   &$112$ &$>250$ \\
$\tilde{P}_U$                &$38$ &$68$ &$81$   &$92$ &$>250$ \\
$\tilde{P}_{SSOR}$        &$38$ &$68$ &$108$   &$82$ &$114$ \\
$\tilde{P}_{ILU}$           &$38$ &$68$ &$83$   &$82$ &$74$ \\
\bottomrule
\end{tabular}\caption{The number of iterations of the preconditioned GMRES method 
  for solving the coupled FSI system with 
  varying fluid density $\rho_f\in\{1.1, 0.11, 0.011, 0.0011, 0.00011\}$ g/cm$^3$.}
\label{tab:unrob}
\end{table}
As observed, the preconditioned Krylov subspace methods do not show 
the robustness with respect to the varying fluid density, i.e., 
the iteration numbers grow more or less when the fluid density decreases. 
In addition, the iteration numbers in the first column 
%conform 
correspond 
to the 
numerical results in \cite{UY14}, where a similar fluid density has been adopted. 
Note that, here we stop the linear 
solver when the error in the GMRES iteration is reduced 
by a factor of $10^{10}$. The numerical results are shown for the first Newton iteration. 

Although we are able to cure the mesh dependence issue 
by using the fully coupled 
monolithic 
AMG methods, 
%when 
provided 
we 
have designed effective smoothers and coarsening strategies for 
such multifield problems, this task in general turns out to 
be nontrivial, see, e.g., \cite{NME:NME3001}. 

Motivated by the above observations, in this work, we aim to 
% construct the preconditioned 
%Krylov subspace methods for the monolithic coupled 
%FSI system with a more robust and efficient preconditioner, 
construct a more robust and efficient preconditioner 
in preconditioned Krylov subspace methods for the monolithic coupled 
FSI system, 
that is based on 
the approximation of 
the direct complete $LDU$ factorization of the coupled system matrix.

The remainder 
of paper 
is organized 
%in the following. 
as follows. 
In Section \ref{sec:model}, 
we set up a FSI model problem for testing our methods. 
Section \ref{sec:lsm} deals with the construction of the robust and efficient 
preconditioner in 
%the 
Krylov subspace 
%method
methods 
for the linearized and discretized 
model problem. Numerical studies are presented in Section \ref{sec:num}. Finally, 
some conclusions are drawn in Section \ref{sec:con}.

%%%%%%%%%%%%%%%%%%%%%%%%%%%%%%%%%%%%%%%%%%%%%%%%%%%%%%%%%%%%%%%%%%%%%%
\section{A model 
%problem setting for the FSI
FSI problem and its discretization
}\label{sec:model}
%Here we set up a model problem for the FSI simulation
\subsection{The geometrical configuration, mappings and kinematics}
\begin{figure}[htbp]
  \centering
  \scalebox{0.3}{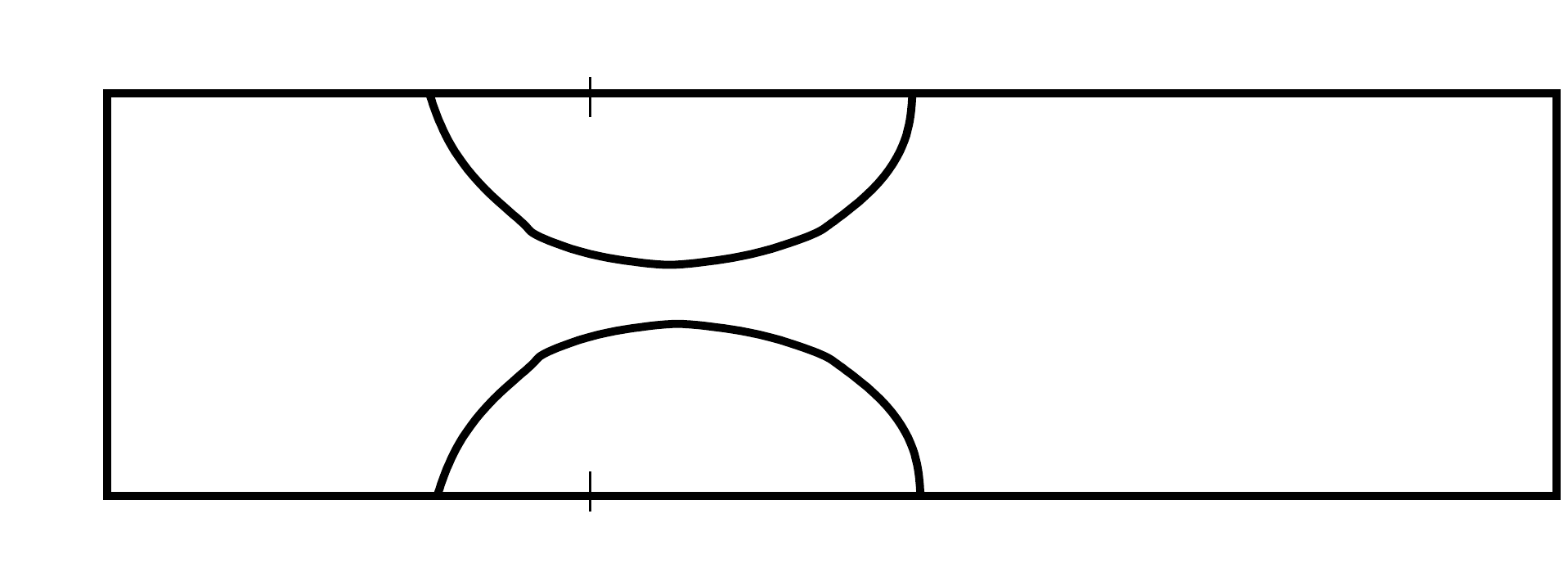}\quad\;
  \scalebox{0.3}{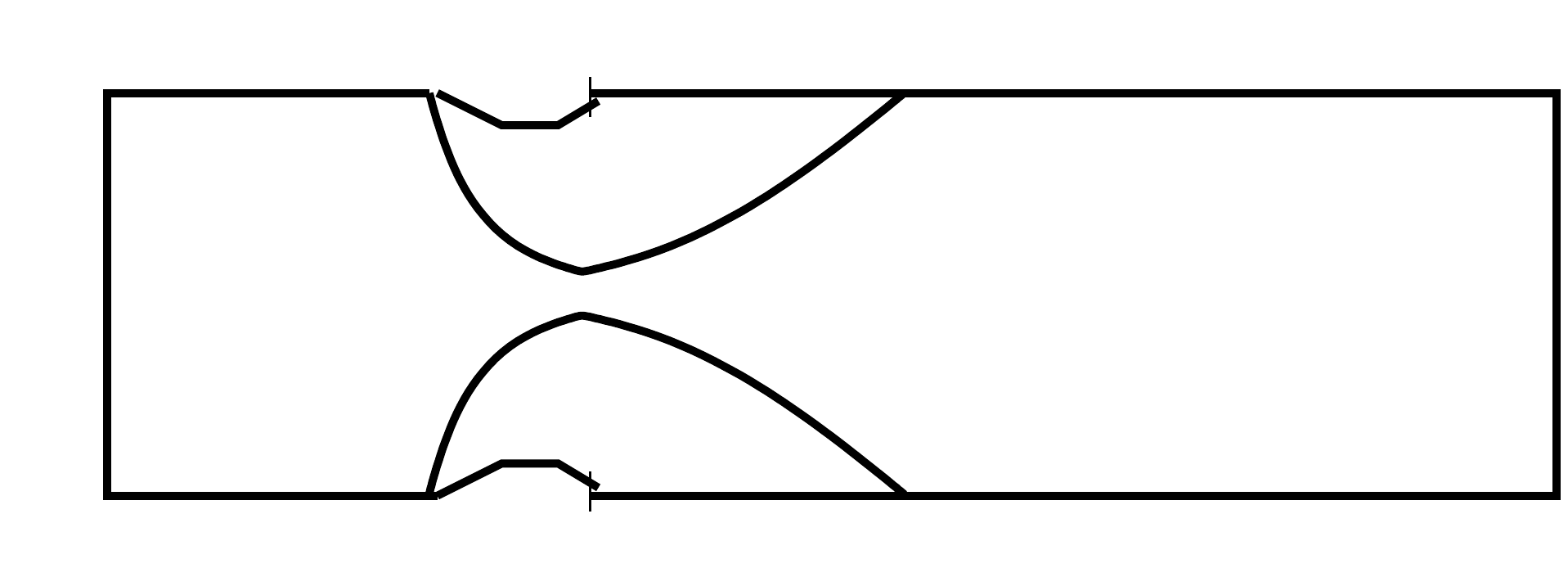}
  \caption{An illustration of the computational FSI domain at the initial time $t=0$ (left) and the current time $t$ (right).}
  \label{fig:dom}
\end{figure}

We consider the FSI domain $\Omega^t$ as a union of the deformable 
fluid domain $\Omega_f^t$ and structure domain $\Omega_s^t$ at the time $t$: 
%$\overline{\Omega^t}:=\overline{\Omega_f^t}\cup\overline{\Omega_s^t}$ 
$\overline{\Omega}^t:=\overline{\Omega}_f^t\cup\overline{\Omega}_s^t$ 
and $\Omega_f^t\cap\Omega_s^t=\emptyset$. The boundary 
of the fluid domain $\partial\Omega_f^t$ 
is decomposed into several parts: 
%$\partial\Omega_f^t:=\overline{\Gamma_{in}}\cup\overline{\Gamma_{out}}\cup\overline{\Gamma_{wall}}\cup\overline{\Gamma^t}$ 
$\partial\Omega_f^t:=\overline{\Gamma}_{in}\cup\overline{\Gamma}_{out}\cup\overline{\Gamma}_{wall}\cup\overline{\Gamma}^t$ 
and $\Gamma_{in}\cap\Gamma_{out}=\Gamma_{in}\cap\Gamma_{wall}=\Gamma_{in}\cap\Gamma^t=\Gamma_{out}\cap\Gamma_{wall}=\Gamma_{out}\cap\Gamma^t=\emptyset$. 
In an analogous way, the boundary $\partial\Omega_f^s$ of the structure domain is decomposed into the following parts: 
%$\partial\Omega_s^t:=\overline{\Gamma_{d}}\cup\overline{\Gamma_{n}^t}\cup\overline{\Gamma^t}$ and $\Gamma_{d}\cap\Gamma_{n}^t=\Gamma_{d}\cap\Gamma^t=\Gamma_{n}^t\cap\Gamma^t=\emptyset$. 
$\partial\Omega_s^t:=\overline{\Gamma}_{d}\cup\overline{\Gamma}_{n}^t\cup\overline{\Gamma}^t$ and $\Gamma_{d}\cap\Gamma_{n}^t=\Gamma_{d}\cap\Gamma^t=\Gamma_{n}^t\cap\Gamma^t=\emptyset$. 
The interface $\Gamma^t$ is defined as the intersection of 
the fluid and structure boundary: $\Gamma^t:=\partial\Omega_f^t\cap\partial\Omega_s^t$. 
At the time level $t=0$, we have the initial (reference) configuration. 
See a schematic illustration in Fig. \ref{fig:dom}.

As usual, we use the Arbitrary-Lagrangian-Eulerian (ALE) mapping 
defined as ${\mathcal A}^t:={\mathcal A}^t(x)= x + d_f(x, t), \;\forall x\in\Omega_f^0$, to 
track the movement of the fluid domain $\Omega_f^0$, where the $d_f:=d_f(x, t),\;\forall x\in\Omega_f^0$ denotes 
the fluid domain displacement; see, e.g., \cite{TH:81, LF99:00, JD04:00}. 
For the structure sub-problem, the Lagrangian mapping ${\mathcal L}^t:={\mathcal L}^t(x)= x + d_s(x, t), \;\forall x\in\Omega_s^0$ is used to track the structure body 
movement; see, e.g., \cite{JB08:00, GAH00:00}. 
In addition, the fluid 
velocity $u:=u(x, t),\;\forall x\in\Omega_f^0$ and 
pressure $p:=p(x, t),\;\forall x\in\Omega_f^0$ are defined via the ALE mapping: 
$u(x, t)=\tilde{u}(\tilde{x}, t)=\tilde{u}({\mathcal A}_f^t(x), t)$, 
$p(x, t)=\tilde{p}(\tilde{x}, t)=\tilde{p}({\mathcal A}_f^t(x), t),\;\forall x\in\Omega_f^0$ and $\tilde{x}={\mathcal A}_f^t(x)\in\Omega_f^t$, where 
$\tilde{u}(\cdot, \cdot)$ and $\tilde{p}(\cdot, \cdot)$ denote 
the fluid velocity and pressure variables under the Eulerian framework; 
see, e.g., \cite{Wick11:00}.

Since we formulate the coupled FSI system on the reference configuration for 
the fluid sub-problem by the ALE mapping 
and structure sub-problem by the Lagrangian mapping, 
we need to introduce the necessary notations 
in the kinematics as used in, e.g., \cite{JB08:00, GAH00:00}. For this, we define 
the fluid and structure deformation gradient tensor by 
$F_f:=F_f(x)=\partial{\mathcal A}^t/{\partial x}=I+\nabla d_f,\;\forall x\in\Omega^0_f$ 
and $F_s:=F_s(x)=\partial{\mathcal L}^t/{\partial x}=I+\nabla d_s,\;\forall x\in\Omega^0_s$, 
respectively. Their determinants are given by
 $J_f=\text{det}F_f$ and $J_s=\text{det}F_s$, accordingly. 

\subsection{A monolithic FSI system on the reference configuration}
After the above preliminary, we formulate the the coupled FSI system in strong form on the 
reference domain: Find $(d_f, u, p, d_s)$ 
such that 
 \begin{subequations}\label{eq:fsi}
    \begin{eqnarray}
      -\Delta d_f=0&{\textup{ in }} \Omega_f^0, \\ %[0.2cm]
      d_f=d_s&{\textup{ on }} \Gamma^0,\\ %[0.2cm]
      \rho_fJ_f\partial_t u+\rho_fJ_f((u- w_f)\cdot F_f^{-1}\nabla) u & \nonumber \\ 
      -\nabla\cdot (J_f\sigma_f(u, p)F_f^{-T})=0 &{\textup{ in }}\Omega_f^0,\\
      \nabla\cdot (\rho_fJ_fF_f^{-1}u)=0&{\textup{ in }}\Omega_f^0, \\ %[0.2cm]
      \rho_s\partial_{tt}d_s-\nabla\cdot (F_sS)=0&{\textup{ in }}\Omega_s^0,\\  
      %-(J_s-1)-(1/\kappa)p_s=0&{\textup{ in }}\Omega_s^0, \\%[0.2cm]
      u=\partial_t d_s&{\textup{ on }}\Gamma^0, \\
      J_f\sigma_{f}(u, p)F_f^{-T}n_f+F_sSn_s=0&{\textup{ on }}\Gamma^0.
    \end{eqnarray}
  \end{subequations}
To complete the system, we prescribe 
the corresponding boundary conditions $d_f=0$ on 
$\Gamma_{in}\cup\Gamma_{wall}\cup\Gamma_{out}$, $u=0$ on $\Gamma_{wall}$, 
$u= g_{in}$ (a given function) on $\Gamma_{in}$ and 
$J_F\sigma_f(u, p)F_f^{-T}n_f = 0$ on $\Gamma_{out}$, 
$d_s=0$ on $\Gamma_d$ and 
$F_sSn_s=0$ on $\Gamma_n$, and the proper initial conditions 
$u(x, 0)=0,\;\forall x\in\Omega_f^0$ and 
$d_s(x,0)=\partial_td_s(x, 0)=0,\;\forall x\in\Omega_s^0$. 

Here, the notations $\rho_f$ and $\rho_s$ 
denote the fluid and structure density, respectively, 
$n_f$ and $n_s$ the fluid and structure outerward unit normal vector, respectively, 
$\sigma_f(u, p):=\mu(\nabla u+\nabla^T)-pI$ the Cauchy stress tensor with the 
dynamic viscosity term $\mu$. 

For the structure, we use the hyperelastic model of the 
%St. Venant 
St.~Venant 
Krichhoff material, for which the second Piola Kirchhoff stress tensor $S$ is 
defined as 
\begin{equation}\label{eq:stvmodel}
  S:=\lambda_s\text{tr}(E_s)I+2\mu_sE_s
\end{equation}
where 
$E_s:=0.5(F_s^TF_s-I)$ denotes the Green-Lagrange strain tensor with 
the Lam\'{e} constant $\lambda_s$ and the shear modulus $\mu_s$.
\begin{remark}
  Note that, in \cite{STurek06, Wick11:00}, 
  a monolithic formulation was introduced and used for the FSI simulation, 
  where the structure 
  velocity variable is introduced to rewrite the structure equation into 
  a system of two first order time dependent equations so that 
  both the fluid and structure sub-problems are rewritten in a monolithic 
  manner on the reference domain $\Omega^0$. In our approach, 
  we keep the form of structure equation, but transform the fluid sub-problem 
  onto the fluid reference domain by the ALE mapping. 
  By this means, we keep the modulus of each sub-problem so that robust solvers 
  for the sub-problem can be directly applied.
  %herein. 
\end{remark}
\subsection{Temporal and spatial discretization and linearization}
Concerning the temporal and spatial discretization, we follow the 
approaches in our previous work \cite{CVULHY1300, UY14}. 
For the time discretization of the fluid and structure sub-problem, we use 
the first order implicit Euler scheme and 
a first order Newmark-$\beta$ scheme, respectively. 
For the spatial discretization of the fluid sub-problem, we use 
the stabilized $P_1-P_1$ finite element discretization with 
standard hat basis functions for both the fluid velocity 
and pressure interpolations. For the mesh movement and structure sub-problem, 
we use the $P_1$ finite element discretization with the 
standard hat function for both the fluid and structure displacement interpolations.
Following the approach in \cite{UY14}, the nonlinearity of the monolithic 
FSI system is handled by Newton's method. 

\section{Monolithic solution methods for the coupled FSI system}\label{sec:lsm}
%In this Section, we will discuss in detail how to construct the robust and 
%efficient preconditioner for the coupled FSI system, based on the 
%complete $LDU$ factorization of the properly reordered system matrix.
\subsection{The modified coupled FSI system}
%Base 
Based 
on our previous work \cite{UY14}, the linearized coupled FSI system 
%to be solved by any linear solver 
is formulated as 
\begin{equation}\label{eq:linsystem1} 
  \left[\begin{array}{cc|ccc|cc}
      A_m^{ii}&A_m^{i\gamma}& & & & & \\ [0.1cm] 
      & I & & & & -I&  \\ [0.1cm] \hline
      B^i_{fm}& B^\gamma_{fm}& -C_f & B^{i}_{1f} & B^{\gamma}_{1f} & & \\ [0.1cm] 
      A^{ii}_{fm}&A^{i\gamma}_{fm}& B^{i}_{2f}&A^{ii}_f&A^{i \gamma}_f&  & \\ [0.1cm]  
      A^{\gamma i}_{fm}&A^{\gamma\gamma}_{fm}&B^{\gamma}_{2f}
      &A^{\gamma i}_f&A^{\gamma\gamma}_f
      &A^{\gamma\gamma}_s&A^{\gamma i}_s\\ [0.1cm] \hline
      & & & & -I & \frac{1}{\Delta t}I & \\ [0.1cm] 
      & & & &    & A^{i\gamma}_s& A^{ii}_s
    \end{array}
    \right]
  \left[\begin{array}{c}
      \Delta d_m^i\\ [0.1cm]  
      \Delta d_m^{\gamma} \\  [0.1cm] \hline
      \Delta p\\ [0.1cm]  
      \Delta u_f^i  \\   [0.1cm]  
      \Delta u_f^{\gamma}\\ [0.1cm] \hline
      \Delta d_s^{\gamma} \\  [0.1cm] 
      \Delta d_s^{i}
    \end{array}\right]   
  =
  \left[\begin{array}{c}
      r_m^i\\ [0.1cm] 
      r_m^{\gamma} \\ [0.1cm] \hline
      r_{p}\\ [0.1cm] 
      r_f^i \\[0.1cm] 
      r_f^{\gamma}\\[0.1cm]  \hline
      r_s^{\gamma}\\[0.1cm] 
       r_s^{i}
    \end{array}\right],
\end{equation}
where the superscript $\gamma$ indicates the 
(nodal) degrees of freedom (DOF) associated to the variables on the interface, 
$i$ denotes the remaining DOF, 
$\beta\alpha$, $\beta, \alpha\in\{\gamma, i\}$, $\beta\neq\alpha$ means the coupling 
between the corresponding interface variables and the remaining. 
Here the second row corresponds to the interface coupling between the mesh movement 
and the structure displacement, the fifth and sixth rows correspond to the 
equivalence of surface tractions and the no-slip interface condition 
from the fluid and structure side, respectively. It is easy to see, on the main diagonal, 
we have stiffness matrices for the mesh movement, fluid and structure sub-problem, 
respectively, and the off-diagonal matrices denote the the coupling among them. 

In order to derive the robust preconditioner for the coupled FSI system in a 
convenient manner, we reoder the 
system (\ref{eq:linsystem1}) by changing the rows and columns, and modify 
some of the matrix entries and the right hand side accordingly. 
This way, we obtain the following equivalent linear system of equations: 
\begin{equation}\label{eq:linsystem1sym} 
  \left[\begin{array}{cc|cc|ccc}
      A_m^{ii}&& & A_{ms}^{i\gamma} & &  & \\ [0.1cm]
      & I &  & -I& &&  \\ [0.1cm] \hline
      & &A^{ii}_s  & & \Delta t A^{i\gamma}_{sf}   & & \\ [0.1cm] 
      & & &\frac{1}{\Delta t}I&   -I &  & \\ [0.1cm] \hline
      A^{\gamma i}_{fm}&A^{\gamma\gamma}_{fm}& A^{\gamma i}_s
      &A^{\gamma\gamma}_s &A^{\gamma\gamma}_f
      &A^{\gamma i}_f &B^{\gamma}_{2f}\\ 
      A^{ii}_{fm}&A^{i\gamma}_{fm}& & &A^{i \gamma}_f&  A^{ii}_f &B^{i}_{2f} \\ [0.1cm] 
      B^i_{fm}& B^\gamma_{fm}&  & & B^{\gamma}_{1f} &B^{i}_{1f}  &-C_f  \\ 
    \end{array}
    \right]
  \left[\begin{array}{c}
      \Delta d_m^i\\ [0.1cm]  
      \Delta d_m^{\gamma} \\  [0.1cm] \hline
      \Delta d_s^{i}\\ [0.1cm] 
      \Delta d_s^{\gamma} \\  [0.1cm] \hline
      \Delta u_f^{\gamma}\\ [0.1cm]  
      \Delta u_f^i  \\   [0.1cm]  
      \Delta p
    \end{array}\right]   
  =
  \left[\begin{array}{c}
      \tilde{r}_m^i\\ [0.1cm] 
      r_m^{\gamma} \\ [0.1cm] \hline
      \tilde{r}_s^{i}\\[0.1cm]  
      r_s^{\gamma}\\[0.1cm] \hline
      r_f^{\gamma}\\ [0.1cm] 
      r_f^i \\[0.1cm] 
      r_{p}
    \end{array}\right], 
\end{equation}
where 
\begin{equation}\label{eq:mrhs}
  \tilde{r}_m^i=r_m^i-A_m^{i\gamma}r_m^\gamma,\;
  \tilde{r}_s^i= r_s^{i} - \Delta t A_s^{i\gamma} r_s^\gamma.
\end{equation}
By this means, we keep the system matrix for the sub-problems on the diagonal 
as symmetric as possible to make our linear solver for the sub-problem 
more efficient, e.g., by 
applying conjugate gradient method with AMG 
preconditioner (see \cite{UH:02}). 
Nevertheless, keeping the symmetry of the sub-problem 
is not mandatory. For instance, for nonsymmetric positive 
systems, a class of AMG methods with special transfers 
base on Schur complements and Galerkin projections are proposed 
in the recent work \cite{NLA:NLA1889}. 
On the other hand, we face the unsymmtry for the fluid sub-problem any way, 
mainly due to the convection and stabilization 
%term. 
terms. 
This issue is 
handled by a class of 
coupled AMG methods \cite{WM04:00,WM06:00,Yang11:00,YH11:00, CVULHY1300}.

For convenience of the following presentation, we rewrite the reordered 
system (\ref{eq:linsystem1sym}) in the following compact form:
\begin{equation}\label{eq:fsieqcomp}
  Kx=b,
\end{equation}
where
\begin{equation}\label{eq:k}
  K= 
  \left[\begin{array}{ccc}
      A_m& A_{ms} & 0\\
      0 &A_{s}& A_{sf} \\
      A_{fm}& A_{fs} & A_{f}\\
    \end{array}
    \right],\;
  x=\left[\begin{array}{c}
      x_m\\
      x_s \\
      x_f\\
    \end{array}
    \right],
  b=\left[\begin{array}{c}
      b_m\\
      b_s \\
      b_f\\
    \end{array}
    \right].
\end{equation}
Here the block matrices and vectors are assigned according to the 
sub-division of the FSI system (\ref{eq:linsystem1sym}).
\subsection{Construct the Schur complement approximation}
To construct the Schur complement approximation of the 
FSI system, we start to perform a $LDU$ factorization for the $3\times 3$ block matrix $K$. 
In this case, the factorization is formulated as 
\begin{equation}\label{eq:fsiredfact3}
  \begin{aligned}
  K=&LDU\\
  :=&\left[\begin{array}{ccc}
      I & 0 & 0\\
      0 & I & 0\\
      A_{fm}A_m^{-1} &  \tilde{A}_{fs}A_s^{-1}  & I  \\
    \end{array}
    \right]
  \left[\begin{array}{ccc}
      A_m & 0 & 0\\
      0 &  A_s  & 0  \\
      0 & 0 & S
    \end{array}
    \right]  
  \left[\begin{array}{ccc}
      I & A_m^{-1}A_{ms} & 0\\
      0 &  I & A_s^{-1}A_{sf} \\
      0 & 0 & I
    \end{array}
    \right]\\
  =&\left[\begin{array}{ccc}
      A_m & 0 & 0\\
      0 & A_s & 0\\
      A_{fm} &  \tilde{A}_{fs}  & S  \\
    \end{array}
    \right]
  \left[\begin{array}{ccc}
      I & A_m^{-1}A_{ms} & 0\\
      0 &  I & A_s^{-1}A_{sf} \\
      0 & 0 & I
    \end{array}
    \right],
  \end{aligned}
\end{equation}
where the fluid Schur complement is formulated as 
\begin{equation}\label{eq:fsiredfactschur}
  S=A_f- \tilde{A}_{fs}A_s^{-1}A_{sf}
\end{equation}
with 
\begin{equation}\label{eq:fsiredfactfs}
  \tilde{A}_{fs}=A_{fs}-A_{fm}A_m^{-1}A_{ms}.
\end{equation}
Inspired by this observation, we propose the following 
FSI preconditioner
\begin{equation}\label{eq:fsiredfactpre}
  \begin{aligned}
  \hat{K}
  =&\hat{L}\hat{D}\hat{U}\\
  :=&\left[\begin{array}{ccc}
      I & 0 & 0\\
      0 & I & 0\\
      A_{fm}\hat{A}_m^{-1} &  \hat{\tilde{A}}_{fs}\hat{A}_s^{-1}  & I  \\
    \end{array}
    \right]
  \left[\begin{array}{ccc}
      \hat{A}_m & 0 & 0\\
      0 &  \hat{A}_s  & 0  \\
      0 & 0 & \hat{S}
    \end{array}
    \right]  
  \left[\begin{array}{ccc}
      I & \hat{A}_m^{-1}A_{ms} & 0\\
      0 &  I & \hat{A}_s^{-1}A_{sf} \\
      0 & 0 & I
    \end{array}
    \right]\\
  =&\left[\begin{array}{ccc}
      \hat{A}_m & 0 & 0\\
      0 & \hat{A}_s & 0\\
      A_{fm} &  \hat{\tilde{A}}_{fs}  & \hat{S}  \\
    \end{array}
    \right]
  \left[\begin{array}{ccc}
      I & \hat{A}_m^{-1}A_{ms} & 0\\
      0 &  I & \hat{A}_s^{-1}A_{sf} \\
      0 & 0 & I
    \end{array}
    \right].
  \end{aligned}
\end{equation}
Here the approximation of the fluid Schur complement is defined as  
\begin{equation}\label{eq:fsiredfactschurapp}
  \hat{S}=A_f- \hat{\tilde{A}}_{fs}\hat{\hat{A}}_s^{-1}A_{sf}
\end{equation}
with 
\begin{equation}\label{eq:fsiredfactfs}
  \hat{\tilde{A}}_{fs}=A_{fs}-A_{fm}\hat{\hat{A}}_m^{-1}A_{ms},
\end{equation}
where 
\begin{equation}\label{eq:fsiredfactasam}
  \hat{\hat{A}}_m = \text{diag}[A_m], \;
      \hat{\hat{A}}_s = \text{diag}[A_s].
\end{equation}
The notation "diag" means the block diagonal of the corresponding matrix from 
the mesh movement and the structure sub-problem, respectively. By 
this means, we are able to construct the fluid Schur complement in an 
explicit way, that corresponds to a full fluid matrix perturbed 
by the matrix from the multiplication of the 
approximated coupling matrices of the fluid and mesh movement 
sub-problem, the fluid and structure sub-problem, and the structure and 
mesh movement sub-problem, respectively. 
\begin{remark}
  In principle, 
  %one 
  we 
  can choose different approximations 
  %for 
  to 
  $A_m$ appearing in the $\hat{L}$, $\hat{D}$, and $\hat{U}$ 
  of (\ref{eq:fsiredfactpre}). 
  In our case, we utilize the same approximation $\hat{A}_m$, i.e., 
  one corresponding AMG iteration is applied to the mesh movement sub-problem 
  for the inverse approximation. The same applies to the 
  approximation for $A_s$. 
\end{remark}
\begin{remark}
  To approximate $A_m$ and $A_s$ appearing in the 
  fluid Schur complement (\ref{eq:fsiredfactschur}), we 
  employ the diagonal of the matrix as an approximation, 
  that turns out to be a rather robust and meanwhile cheap 
  approximation in our applications. However, in principle, the methodology here can 
  be extended to other situations, where the approximation 
  for the inverse of $A_m$ and $A_s$ can be computed 
  explicitly in another cheap way. 
\end{remark}
\begin{remark}
  The constructure of the fluid Schur complement approximation (\ref{eq:fsiredfactschurapp}) 
turns out to be a fairly cheap operation. Since the matrices 
$A_{fs}$, $A_{ms}$ and $A_{sf}$ have very sparse non-zero 
pattern corresponding to the coupling conditions among the 
mesh movement, fluid and structure sub-problems on the interface only, 
the cost of multiplication between the matrices, and between 
the matrix and the diagonal 
of the matrix is rather cheap. In addition, the sparse 
matrix $A_{fm}$ couples the mesh movement and the fluid sub-problem 
in the fluid reference 
%domain, so 
domain. Thus, 
it has more entries than the interface 
coupling matrices. Finally, we only need to modify the entries in the $A_f$ 
and $B_{1f}$ blocks of the fluid matrix in (\ref{eq:linsystem1sym}) 
in order to construct the fluid Schur complement. 
By neglecting the perturbation $-A_{fm}A_m^{-1}A_{ms}$ in (\ref{eq:fsiredfactschur}), 
a cheaper inexact approximation for the Schur complement $S$ is 
obtained. However, this approximation turns out to 
be 
too rough to get 
the robustness of the preconditioner. Thus it will not be discussed in the following. 
\end{remark}
\begin{remark}
  There are two other direct ways to construct the Schur complement 
  for the coupled FSI system. 

  One way is to construct the structure 
  Schur complement:
  \begin{equation}\label{eq:fsischurprestruc}
    S = A_s-A_{sf}A_{f}^{-1}\tilde{A}_{fs}=A_s-A_{sf}A_{f}^{-1}(A_{fs}-A_{fm}A_m^{-1}A_{ms}),
  \end{equation}
  based on the following $LDU$ factorization of the original coupled 
  system:
  \begin{equation*}
    \begin{aligned}
      &\left[\begin{array}{ccc}
          A_m& 0 & A_{ms}\\
          A_{fm}& A_{f}& A_{fs} \\
          0 &A_{sf}& A_s
        \end{array}
        \right]\\
      =& 
      \left[\begin{array}{ccc}
          I & 0 & 0\\
          A_{fm}A_m^{-1}&  I & 0  \\
          0 &A_{sf}A_f^{-1}& I
        \end{array}
        \right]
      \left[\begin{array}{ccc}
          A_m & 0 & 0\\
          0 &  A_f  & 0  \\
          0 & 0 & S
        \end{array}
        \right]  
      \left[\begin{array}{ccc}
          I & 0 & A_m^{-1}A_{ms}\\
          0 &  I & A_f^{-1}\tilde{A}_{fs} \\
          0 & 0 & I
        \end{array}
        \right],
    \end{aligned}
  \end{equation*}
  where $\tilde{A}_{fs}=A_{fs}-A_{fm}A_m^{-1}A_{ms}$. 
  
  The other way is to construct the mesh movement Schur complement:
  \begin{equation}\label{eq:fsischurpremesh}
    S=A_m-A_{mf}A_{fm}=A_m+A_{ms}A_s^{-1}A_{sf}(A_f -A_{fs}A_s^{-1}A_{sf})^{-1}A_{fm},
  \end{equation}
  based on the following $LDU$ factorization of the reordered FSI system:
  \begin{equation*}
    \begin{aligned}
      &\left[\begin{array}{ccc}
          A_s& A_{sf} & 0\\
          A_{fs} &A_{f}& A_{fm} \\
          A_{ms}& 0 & A_{m}\\
        \end{array}
        \right]\\
      =&\left[\begin{array}{ccc}
          I & 0 & 0\\
          A_{fs}A_s^{-1} & I & 0\\
          A_{ms}A_s^{-1} &  A_{mf} & I  \\
        \end{array}
        \right]
      \left[\begin{array}{ccc}
          A_s & 0 & 0\\
          0 &  \tilde{A}_f & 0  \\
          0 & 0 & S
        \end{array}
        \right]
      \left[\begin{array}{ccc}
          I & A_s^{-1}A_{sf} & 0\\
          0 &  I & \tilde{A}_f^{-1}A_{fm} \\
          0 & 0 & I
        \end{array}
        \right],
    \end{aligned}
  \end{equation*}
where $\tilde{A}_f=A_f  - A_{fs}A_s^{-1}A_{sf}$ and 
$A_{mf}= -{A}_{ms}A_s^{-1}A_{sf}\tilde{A}_f^{-1}$. 

However, none of these two Schur complements (\ref{eq:fsischurprestruc}) 
and (\ref{eq:fsischurpremesh}) is cheaper than 
the fluid Schur complement (\ref{eq:fsiredfactschur}) to 
approximate. From now on, 
we only consider the preconditioner $\hat{K}$ constructed in (\ref{eq:fsiredfactpre}).
\end{remark}
\subsection{The preconditioning steps}
The preconditioning operation in the 
preconditioned Krylov subspace methods, e.g., 
the preconditioned GMRES \cite{Saad86} or the flexible preconditioned 
GMRES \cite{Saad1993}, requires the evaluation $x=\hat{K}^{-1}r$ for a given 
vector $r=[r_m^T, r_s^T, r_f^T]^T$. One inverse operation 
contains five steps as indicated in Algorithm \ref{alg:fsipre}.
\begin{algorithm}
   \caption{Evaluation of $x=\hat{K}^{-1}r$}\label{alg:fsipre} 
   \begin{algorithmic}[1]
     \STATE{$\tilde{x}_m=\hat{A}_m^{-1}r_m$,
     }
     \STATE{$\tilde{x}_s=\hat{A}_s^{-1}r_s$,
     }
     \STATE{$x_f=\hat{S}^{-1}(r_f-A_{fm}\tilde{x}_m-\hat{\tilde{A}}_{fs}\tilde{s}_s)$,
     }
     \STATE{$x_s=\tilde{x}_s-\hat{A}^{-1}_sA_{sf}x_f$,
     }
     \STATE{$x_m=\tilde{x}_m-\hat{A}^{-1}_mA_{ms}x_s$.
     }
   \end{algorithmic}
\end{algorithm}
As observed, we need to evaluate the inverse of $\hat{A}_m$ and $\hat{A}_s$ applied to 
the corresponding vectors twice, and the inverse evaluation of the approximated fluid 
Schur complement $\hat{S}$ applied to a vector once, that is the most 
expensive operation in the preconditiong steps. 
All the evaluation is computed by 
applying one 
%iteration 
cycle 
of a special class of AMG methods \cite{FK98:00, WM04:00, YH11:00, Yang11:00} to the corresponding sub-problem with $0$ 
%starting value, 
as initial guess, 
that turns out to 
be sufficient to evaluate the inverse approximation in the preconditioning steps.   
\subsection{Some remarks on other related preconditioners}
%We like to 
In this subsection, we make some remarks on other known preconditioned 
Krylov subspace methods for the coupled FSI system.
%, that are closely related to our approach. 
%For this, we shortly 
%point out the difference of our proposed method in contrast to 
%the other known methods to our knowledge. 
%Since there exist other different solvers we may be not aware of 
%or less related to the content of this work, our comparison is by far not 
%complete. 
For instance, in \cite{EPFL11, ABXCC10}, the 
domain decomposition based additive Schwarz preconditioners have been 
used in the parallel FSI solvers. 
\begin{remark}
  In \cite{Heil20041}, the author considered the preconditioned 
  Krylov subspace method 
  for the FSI problem in 2D using a 1D model of the wall deformation. Starting 
  from the three block-triangular approximations (as preconditioners) 
  of the original linearized system, the "sup", the 
  "sub" and the "diag", the author then used a global 
  pressure Schur complement preconditioner \cite{ST99:00} to replace 
  the Navier-Stokes block, in which the Elman's {\it{BFBt}} 
  approximation \cite{EH99} 
  to the fluid pressure Schur complement is employed in order to 
  reduce the computational cost. Later on, in \cite{Muddle20127315}, 
  the authors considered a FSI preconditioner for a 
  $4\times 4$ linearized FSI system by replacing the bottom-right 
  $2\times 2$ block in the Jacobian with the so-called 
  pseudo-solid preonditioner. 
  In our approach, we construct the preconditioner based the 
  complete $LDU$ factorization for the the linearized $3\times 3$ coupled 
  system and the approximated Schur complement 
  itself is on the sub-problem level, corresponding to the perturbed fluid 
  sub-problem, see (\ref{eq:fsiredfactschurapp}), for which we have 
  very efficient and robust AMG method 
  \cite{WM04:00, WM06:00, Yang11:00, CVULHY1300} to perform 
  the inverse operation. 
\end{remark}
\begin{remark}
  In \cite{Badia20084216}, the authors considered 
  the domain decomposition 
  Dirichlet-Neumann, the ILU and the inexact 
  block-LU factorization based preconditioners 
  for the FSI system linearized by 
  the fixed point algorithm for both the domain dependence 
  and the convective term. Thus the linear system therein has 
  slightly different structure as we consider in this work. Herein, 
  more coupling matrix blocks come into the the system 
  due to the nonlinearity of the domain movement and the 
  convective term treated in an all-at-once manner 
  using Newton's method. In their inexact block-LU factorization 
  based preconditioner, the inverse of the fluid matrix appears several 
  times in the block structure, that is approximated 
  by Neumann expansion technique. The similar technique was also used in 
  \cite{BS08}. In our approach, 
  by reordering the coupled system we arrive at the fluid 
  Schur complement, that involves the inverse 
  of the perturbed fluid matrix only once. Furthermore, one 
  complete preconditioning step involves the 
  inverse of the mesh movement and structure matrices twice, that are in general 
  cheaper to approximate than the inverse of the fluid sub-problem. 
\end{remark}
\begin{remark}
  In \cite{NME:NME3001}, the block Gauss-Seidel preconditioned Krylov 
  subspace and the fully coupled FSI AMG methods are proposed, 
  which are based on the smoothed aggregation multigrid method for each sub-problem, 
  employed either in the preconditioning step or in the smoothing step. 
\end{remark}
%To sum, 
Summarizing, 
compared to the others, we propose a preconditioner for the linearized 
FSI system, 
based on the complete $LDU$ factorization of the reordered $3\times 3$ system block matrix. 
The Schur complement itself corresponds to a fluid sub-problem, perturbed 
by a sparse matrix coming from the multiplication of 
the corresponding coupling matrices from the mesh movement, fluid and structure sub-problems, that is 
approximated and constructed in an explicit way. 
%\section{Hierarchical matrix approximation for finite element discretiation}
%In order to further reduce the computaional cost, we will consider the hierachical matrix 
%technique to obtain the inverse of the matrix with sparse approximation. 

\section{Numerical experiments}\label{sec:num}
In this Section, we 
%like 
would like 
to demonstrate the robustness of the 
preconditioner in the preconditioned GMRES and 
flexible preconditioned GMRES methods 
for solving the FSI problem. For this reason, we test the algorithm 
for the FSI problem on three consecutively refined finite element meshes. 
We compare the iteration numbers of the preconditioned Krylov subspace methods, 
with varying mesh size, fluid density, and time step size. 
\subsection{Geometry, meshes, material parameters and boundary conditions}
\begin{figure}[htbp]
  \centering
  \includegraphics[scale=0.4]{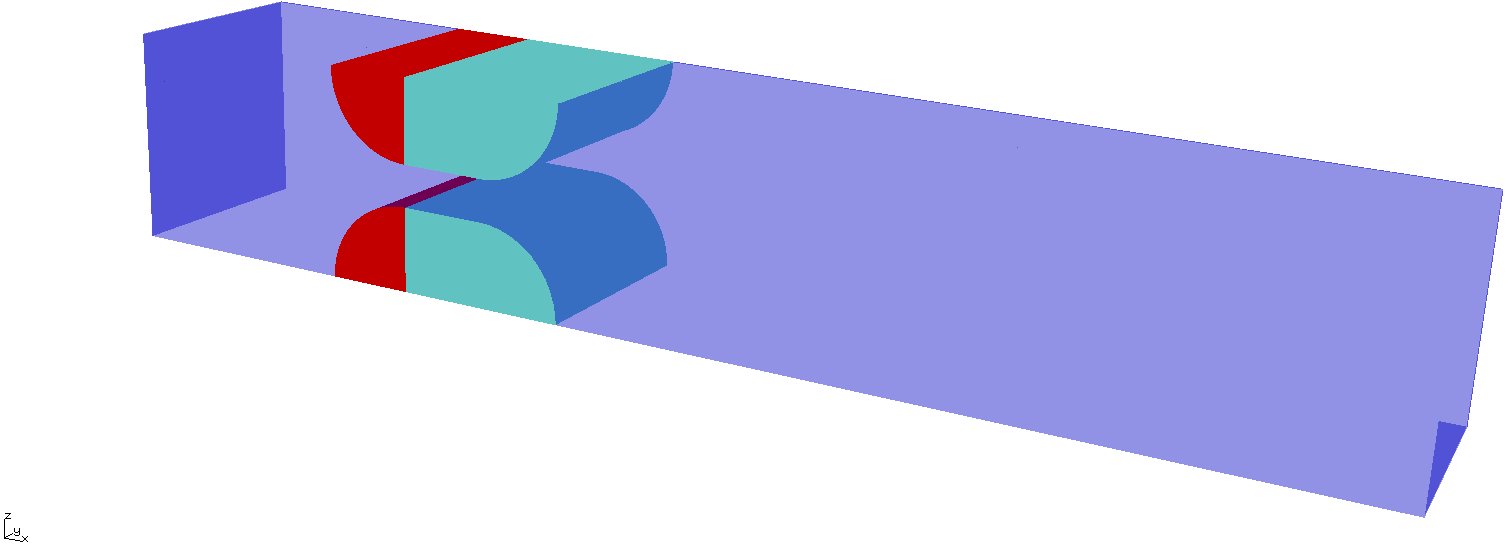}
  \caption{The configuration of the geometry.}
  \label{fig:geo_model}
\end{figure}

\begin{figure}[htbp]
  \centering
  \includegraphics[scale=0.3]{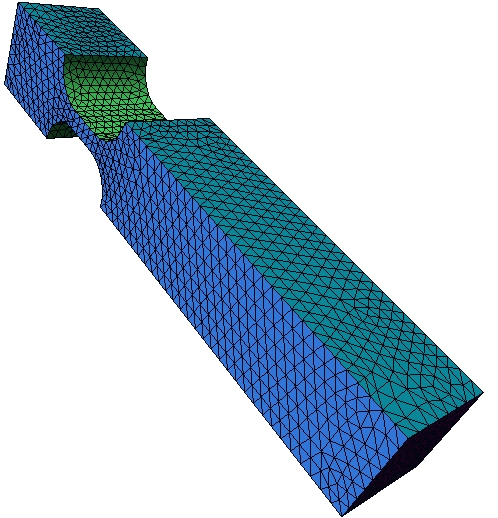}
  \includegraphics[scale=0.3]{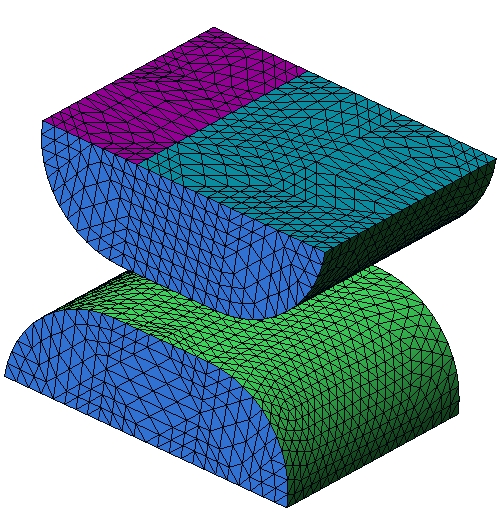}
  \caption{The mesh generated for the fluid domain (left) and the structure domain (right).}
  \label{fig:geo_mesh}
\end{figure}

First we 
%demonstrate 
describe 
the geometry for the 
FSI simulation in Fig. \ref{fig:geo_model}. 
The channel has an obstacle inside, where the 
$x$-, $y$- and $z$-coordinates represent the lateral, 
anterior-posterior, and the vertical directions. respectively. The channel 
has the size $[0, 12]$ cm, $[0, 2]$ cm and $[0, 2]$ cm, 
in the $x$-, $y$- and $z$-direction, respectively. The obstacle is 
composed of four quarter cylinders with radius $0.8$ cm, and two cubes inserted 
in between. The FSI 
interaction occurs on the obstacle surface inside the channel, when the flow 
goes from the left to right in the lateral direction. The finite 
element meshes are generated using Netgen \cite{JS97:00}, where the conforming 
grids on the interface are 
%guaranteed. See 
guaranteed, see 
a mesh example in Fig.\ref{fig:geo_mesh}. 
The information concerning the number of nodes (\#Nod), tetrahedral elements (\#Tet), 
and degrees of freedom (\#Dof) on the coarse mesh (C), intermediate mesh (I) and 
fine mesh (F) is summarized in Table \ref{tab:mesh_info}.
\begin{table}[ht!]
\centering
\begin{tabular}{crrr}
\toprule
&  \#Nod &  \#Tet    & \#Dof \\
\midrule
Coarse mesh (C) & $837$    & $2679$  & $5131$ \\
Intermediate mesh (I) & $5002$ & $21432$  &$31178$  \\
Fine mesh (F) & $34041$  & $171456$  &$214223$ \\
\bottomrule\end{tabular}\caption{Three finite element meshes.}
\label{tab:mesh_info}
\end{table}

We use the nonlinear isotropic and homogeneous hyperelastic model of the 
%St. Venant 
St.~Venant 
Krichhoff material, where the elastic constants in (\ref{eq:stvmodel}) 
are $\lambda_s=1.73$e+$06$ dyne/cm$^2$ 
and $\mu_s=1.15$e+$06$ dyne/cm$^2$. The density of the 
structure is $\rho_s=1$ g/cm$^3$. The fluid kinematic viscosity is $\nu=0.1568$ cm$^2$/s. 
The fluid density is $\rho_f\in\{1.1, 0.11, ..., 0.00011\}$ g/cm$^3$ for our 
testing purpose, that covers a large range of flows, e.g., the water, blood and air flow. 
Note that, for numerical studying purpose, we successively decrease the 
fluid densities by factor of $10$. 

The structure is fixed on the boundaries except that 
part of the top and bottom in the vertical directions 
is assigned a homogeneous Neumann boundary condition as indicated by the purple 
color in Fig.\ref{fig:geo_mesh}. The inflow boundary condition is $u= 30$ cm/s on $\Gamma_{in}$.
For the outflow on $\Gamma_{out}$, 
we use "doing-nothing" condition, i.e., homogeneous Neumann boundary condition. 
On $\Gamma_{wall}$, we use homogeneous Dirichlet boundary condition. 
For testing purpose, we use the time step size $\Delta t\in\{1.25, 0.125, ..., 0.00125\}$ ms. 

\subsection{Numerical studies on the robustness with respect to fluid density, 
  mesh size and time step size}
To study the robustness of the preconditioned Krylov subspace methods with respect to 
the varying fluid density, mesh size and time step size, we test the algorithm on 
three meshes (see Table. \ref{tab:mesh_info}), 
with different fluid densities $\rho_f\in\{1.1, 0.11, ..., 0.00011\}$ g/cm$^3$, 
and different time step size $\Delta t\in\{1.25, 0.125, ..., 0.00125\}$ ms. 
We uses the preconditioned GMRES and flexible GMRES as 
the outer acceleration. The iteration numbers of the preconditioned Krylov subspace methods 
are prescribed in Table \ref{tab:t125} -- \ref{tab:t000125}, 
for solving the linearized FSI system at the first Newton iteration on the first time step. For 
other Newton iterations and time steps, we observe very similar results. 
\begin{table}[ht!]
\centering
\begin{tabular}{cccccccccccccccc}
\toprule
&\multicolumn{15}{c}{\#it}  \\
 \cmidrule(r){2-16} 
 $\rho_f$
 &\multicolumn{3}{c} {1.1} &\multicolumn{3}{c} {0.11}
&\multicolumn{3}{c} {0.011}&\multicolumn{3}{c} {0.0011}
 &\multicolumn{3}{c} {0.00011} \\ 
\cmidrule(r){2-4} \cmidrule(r){5-7}  \cmidrule(r){8-10}  \cmidrule(r){11-13} 
\cmidrule(r){14-16} 
& C    & I & F    & C    & I & F & C    & I & F & C    & I & F& C    & I & F \\
\midrule
Pre\_GMRES        &$7$ &$13$ & $18$   &$6$ &$11$ & $17$    &$5$ &$10$  & $14$ &$4$ &$8$  & $12$  &$3$ &$6$  & $9$ \\
Pre\_FGMRES     &$5$ &$21$ & $18$   &$6$ & $22$ & $38$   &$6$ &$22$ & $38$ &$6$ &$22$ & $38$ &$6$ &$22$ & $38$ \\
\bottomrule
\end{tabular}\caption{The iteration numbers of the preconditioned GMRES and 
flexible GMRES with the time step size $\Delta t = 1.25$ ms.}
\label{tab:t125}
\end{table}
\begin{table}[ht!]
\centering
\begin{tabular}{cccccccccccccccc}
\toprule
&\multicolumn{15}{c}{\#it}  \\
 \cmidrule(r){2-16} 
 $\rho_f$
 &\multicolumn{3}{c} {1.1} &\multicolumn{3}{c} {0.11}
&\multicolumn{3}{c} {0.011}&\multicolumn{3}{c} {0.0011}
 &\multicolumn{3}{c} {0.00011} \\ 
\cmidrule(r){2-4} \cmidrule(r){5-7}  \cmidrule(r){8-10}  \cmidrule(r){11-13} 
\cmidrule(r){14-16} 
& C    & I & F    & C    & I & F & C    & I & F & C    & I & F& C    & I & F \\
\midrule
Pre\_GMRES        &$3$ &$5$ & $8$   &$3$ &$5$ & $8$    &$3$ &$5$  & $5$ &$2$ &$4$  & $5$  &$2$ &$4$  & $5$ \\
Pre\_FGMRES     &$3$ &$7$ & $8$   &$4$ & $8$ & $12$   &$4$ &$9$ & $14$ &$4$ &$9$ & $15$ &$4$ &$9$ & $15$ \\
\bottomrule
\end{tabular}\caption{The iteration numbers of the preconditioned GMRES and 
flexible GMRES with the time step size $\Delta t = 0.125$ ms.}
\label{tab:t0125}
\end{table}
\begin{table}[ht!]
\centering
\begin{tabular}{cccccccccccccccc}
\toprule
&\multicolumn{15}{c}{\#it}  \\
 \cmidrule(r){2-16} 
 $\rho_f$
 &\multicolumn{3}{c} {1.1} &\multicolumn{3}{c} {0.11}
&\multicolumn{3}{c} {0.011}&\multicolumn{3}{c} {0.0011}
 &\multicolumn{3}{c} {0.00011} \\ 
\cmidrule(r){2-4} \cmidrule(r){5-7}  \cmidrule(r){8-10}  \cmidrule(r){11-13} 
\cmidrule(r){14-16} 
& C    & I & F    & C    & I & F & C    & I & F & C    & I & F& C    & I & F \\
\midrule
Pre\_GMRES        &$4$ &$4$ & $4$   &$4$ &$5$ & $4$    &$4$ &$5$  & $4$ &$3$ &$4$  & $3$  &$3$ &$4$  & $3$ \\
Pre\_FGMRES     &$6$ &$6$ & $5$   &$6$ & $7$ & $6$   &$6$ &$8$ & $7$ &$6$ &$8$ & $7$ &$6$ &$8$ & $7$ \\
\bottomrule
\end{tabular}\caption{The iteration numbers of the preconditioned GMRES and 
flexible GMRES with the time step size $\Delta t = 0.0125$ ms.}
\label{tab:t00125}
\end{table}
\begin{table}[ht!]
\centering
\begin{tabular}{cccccccccccccccc}
\toprule
&\multicolumn{15}{c}{\#it}  \\
 \cmidrule(r){2-16} 
 $\rho_f$
 &\multicolumn{3}{c} {1.1} &\multicolumn{3}{c} {0.11}
&\multicolumn{3}{c} {0.011}&\multicolumn{3}{c} {0.0011}
 &\multicolumn{3}{c} {0.00011} \\ 
\cmidrule(r){2-4} \cmidrule(r){5-7}  \cmidrule(r){8-10}  \cmidrule(r){11-13} 
\cmidrule(r){14-16} 
& C    & I & F    & C    & I & F & C    & I & F & C    & I & F& C    & I & F \\
\midrule
Pre\_GMRES        &$4$ &$5$ & $4$   &$4$ &$5$ & $5$    &$4$ &$5$  & $5$ &$3$ &$4$  & $4$  &$3$ &$4$  & $4$ \\
Pre\_FGMRES     &$5$ &$6$ & $6$   &$6$ & $7$ & $7$   &$6$ &$8$  & $8$ &$6$ &$8$   & $8$  &$6$ &$8$ & $8$ \\
\bottomrule
\end{tabular}\caption{The iteration numbers of the preconditioned GMRES and 
flexible GMRES with the time step size $\Delta t = 0.00125$ ms.}
\label{tab:t000125}
\end{table}

The observations are summarized in the following. 
We observe, that the preconditioned GMRES (Pre\_GMRES) method performances 
better than the preconditioned flexible GMRES (Pre\_FGMRES) 
method for our testing problem, i.e., 
the latter needs more iterations for all the test cases. The gap becomes smaller 
when the time step size is refined. Furthermore, the later shows to be more 
sensitive to the fluid density changes, i.e., when the density decreases, the iteration numbers 
very slightly increase. However, for the Pre\_GMRES method, we need fewer iterations when the 
fluid density decreases. In addition, as well known, when 
$\rho_f={\mathcal O}(\rho_s)$ (in our case, $\rho_f\approx \rho_s=1$ g/cm$^3$), 
the so-called added-mass (see, e.g., \cite{HJPMRO95}) plays 
important effect to the normal partitioned FSI solvers, that are usually become 
slow in such a situation; see \cite{Causin04}. However, our preconditioner shows 
the robustness with respect to such added-mass effect; see the 
iteration numbers corresponding to $\rho_f=1.1$ g/cm$^3$ in 
Table \ref{tab:t0125} -- \ref{tab:t000125}.

In any case, we need much fewer iterations than the results in Table \ref{tab:unrob} using 
the old preconditioners, i.e.. our new preconditioner is much more 
efficient than the old ones with respect to the varying fluid densities. By this 
means, we overcome the robustness issues related to the fluid density 
using the old preconditioners.  

Concerning the mesh dependence, from the results in Table \ref{tab:t125}, we 
see slightly increased iterations with the mesh refinement (C-I-F) 
for both methods. However, from the 
results in Table \ref{tab:t0125} -- \ref{tab:t000125} with reduced time step size, 
we observe almost the same 
iterations on the three mesh levels for each method. Compared to the 
results in \cite{UY14} using the old preconditioners, these new results demonstrate 
the more robustness of our new preconditioner with respect to the mesh size. 

Furthermore, regarding the dependence on the time step size, except the result in 
Table \ref{tab:t125}, where a large time step size $\Delta t =1.25$ ms is used in order to 
test the algorithm, that is usually much larger than required in the real simulation, we observe 
more or less similar iterations in Table \ref{tab:t0125} -- \ref{tab:t000125}, 
for a large range of time step size in each method. So, our new preconditioner shows the 
robustness with respect to the time step size. 

We emphasize that, nearly the same complexity of the preconditioning step 
in the new preconditioner is needed in comparison with the old 
preconditioners $\tilde{P}_{SSOR}$ and $\tilde{P}_{ILU}$. 
However, the iteration numbers of the Krylov subspace methods 
with the new preconditioner are approximately reduced by a factor of $10$, 
see Table \ref{tab:unrob} and Table \ref{tab:t0125}. 
%We observe the computational time is approximately reduced by a factor of $10$. 
\subsection{Visualization of the numerical solutions}
\begin{figure}[htbp]
  \centering
  \includegraphics[scale=0.5]{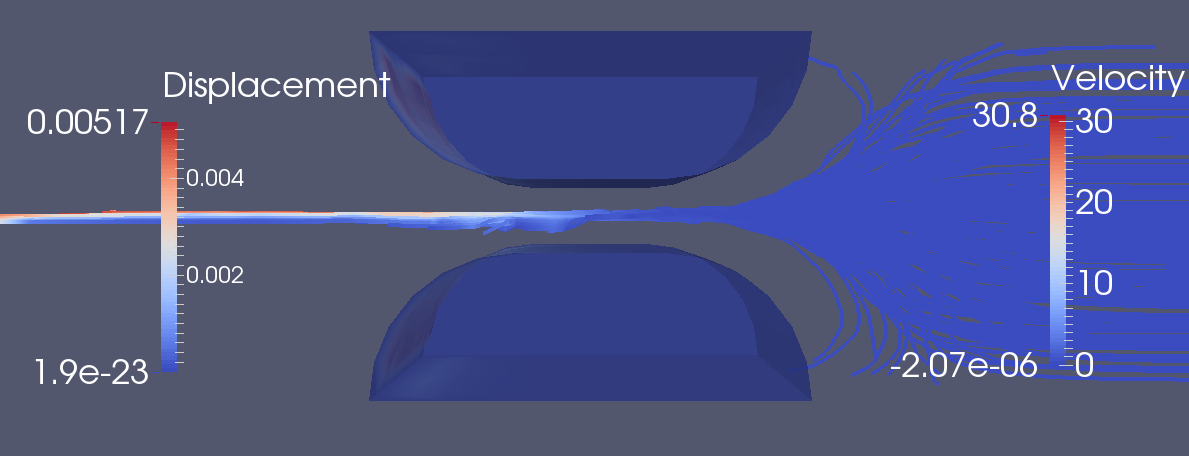}
  \includegraphics[scale=0.5]{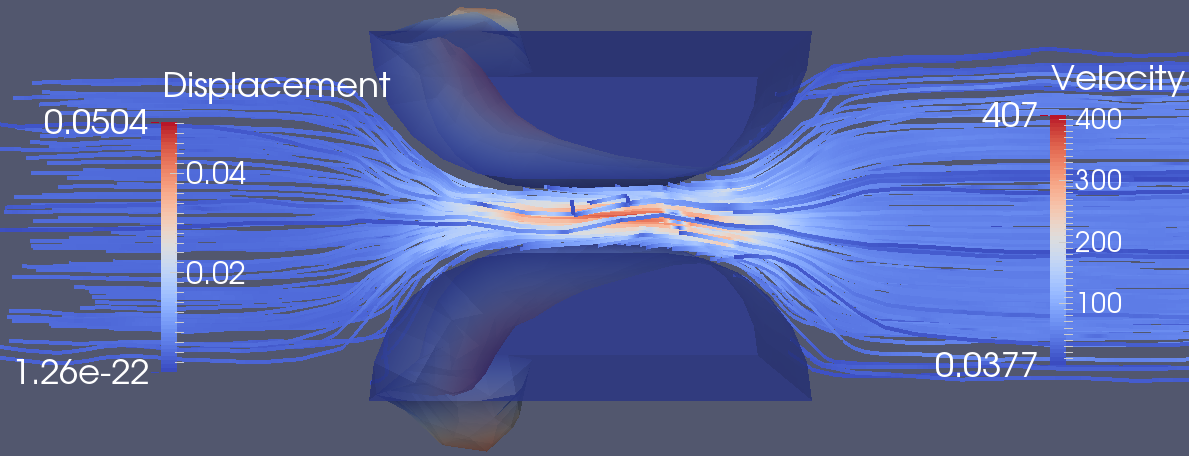}
  \includegraphics[scale=0.5]{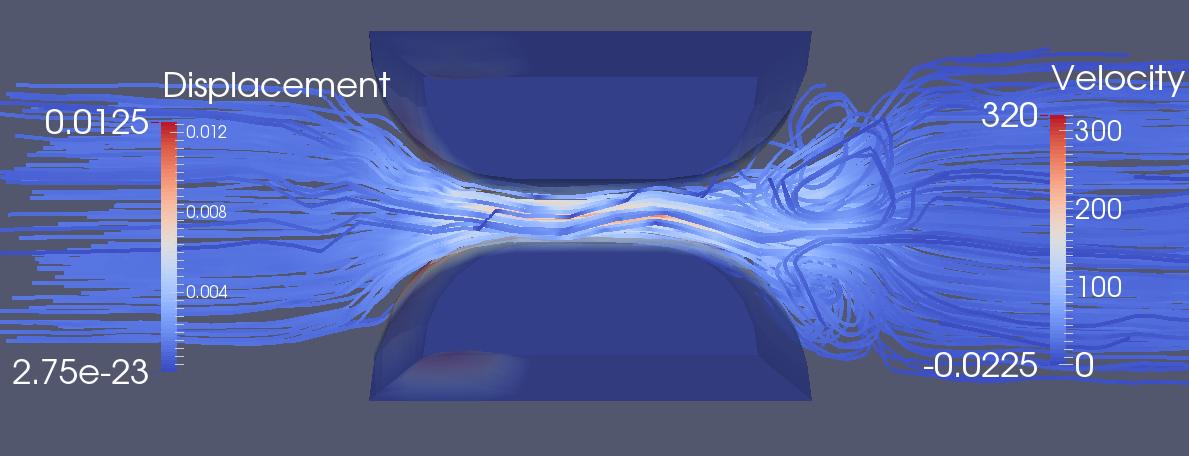}
  \caption{The fluid velocity and structure displacement fields of the FSI simulation with $\rho_f=1.1$ g/cm$^3$ at different time levels: $t=0.125$ ms, $t=7$ ms and $t=14.5$ ms (from top to bottom).}\label{fig:fsi_sol_1}
\end{figure}

\begin{figure}[htbp]
  \centering
  \includegraphics[scale=0.5]{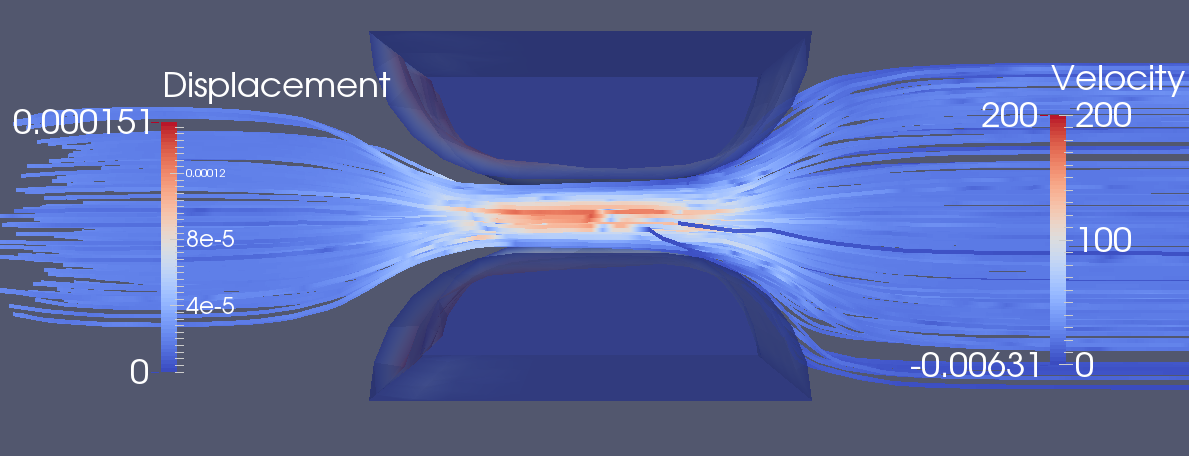}
  \includegraphics[scale=0.5]{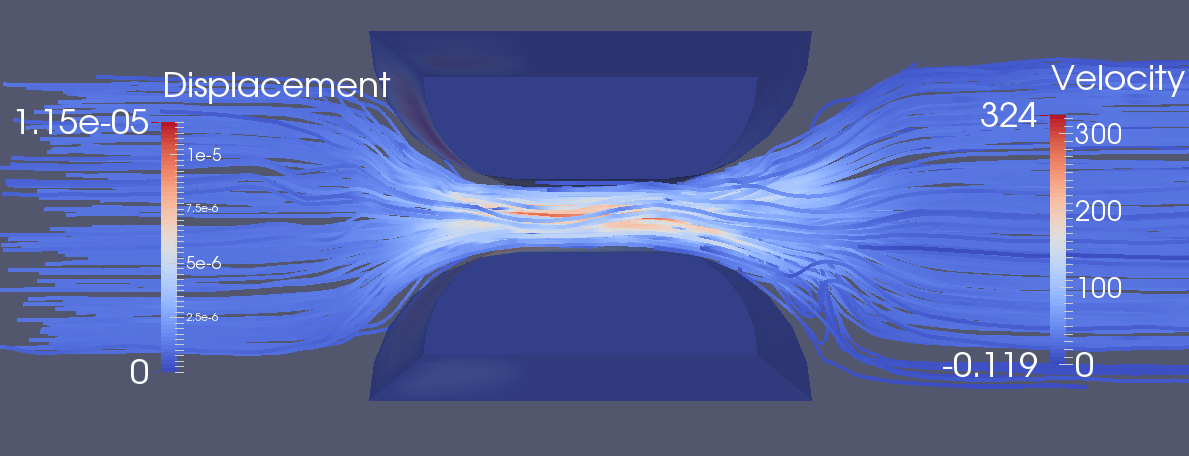}
  \includegraphics[scale=0.5]{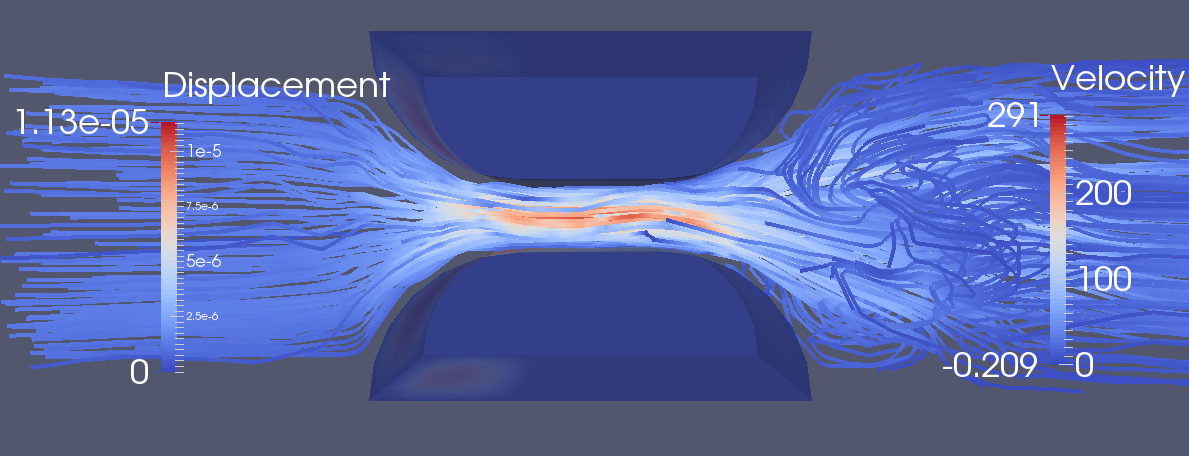}
  \caption{The fluid velocity and structure displacement fields of the FSI simulation with $\rho_f=0.0011$ g/cm$^3$ at different time levels: $t=0.125$ ms, $t=14.5$ and $t=40.5$ ms (from top to bottom).}\label{fig:fsi_sol_2}
\end{figure}

To illustrate the numerical solutions, we 
plot the streamlines of the velocity fields behind the obstacles 
and the structure deformations from the FSI simulation using 
different density $\rho_f=1.1$ g/cm$^3$ (close to water) 
and $\rho_f=0.0011$ g/cm$^3$ (close to air) 
in Fig. \ref{fig:fsi_sol_1} and Fig. \ref{fig:fsi_sol_2}, respectively. The simulation 
results at different time level $t\in\{0.125, 7, 14.5\}$ ms and $t\in\{0.125, 14.5, 40.5\}$ ms 
are shown from top to bottom 
in Fig. \ref{fig:fsi_sol_1} and Fig. \ref{fig:fsi_sol_2}, respectively. 
For visualization purpose, the deformation of the strcuture body is enlarged 
by a factor of $10$. We observe larger structure 
deformation in the FSI simluation with $\rho_f=1.1$ g/cm$^3$ 
than with  $\rho_f=0.0011$ g/cm$^3$. 
In addition, we also obverse some vorticities behind 
the obstacle at late time level, e.g., at $t=14.5$ ms for $\rho_f=1.1$ g/cm$^3$ 
and $t=40.5$ ms for $\rho_f=0.0011$ g/cm$^3$. 
From the numerical experiments, we also observe, that the velocity speed of 
the air flow reaches much higher level than the water flow at the very first 
time step. 

\section{Conclusions}\label{sec:con}
In this work, we develop a new preconditioner for the coupled FSI 
%system, 
sysem of finite element equations 
based on 
the proper approximation of the 
complete $LDU$ factorization of the system matrix. From our numerical studies, 
the preconditioner shows 
the more robustness with respect to the mesh size, varying fluid density, and in addition, 
to the time step size, without resorting to the more involved fully coupled FSI 
monolithic multigrid methods. 
Thus, the new method may reduce the 
implementation and computational complexity significantly. 

\bibliography{RobustFSI}

\newcommand{\noopsort}[1]{} \newcommand{\printfirst}[2]{#1}
  \newcommand{\singleletter}[1]{#1} \newcommand{\switchargs}[2]{#2#1}
\begin{thebibliography}{10}

\bibitem{PSV14}
{\sc A.~Aposporidis, P.~S. Vassilevski, and A.~Veneziani}, {\em Multigrid
  preconditioning of the non-regularized augmented {B}ingham fluid problem},
  Electron. Trans. Numer. Anal., 41 (2014), pp.~42--61.

\bibitem{AO89I}
{\sc O.~Axelsson and P.~S. Vassilevski}, {\em {Algebraic multilevel
  preconditioning methods. I}}, Numer. Math., 56 (1989), pp.~157--177.

\bibitem{AO90}
\leavevmode\vrule height 2pt depth -1.6pt width 23pt, {\em Algebraic multilevel
  preconditioning methods, {II}}, SIAM J. Numer. Anal., 27 (1990),
  pp.~1569--1590.

\bibitem{Badia20084216}
{\sc S.~Badia, A.~Quaini, and A.~Quarteroni}, {\em Modular vs. non-modular
  preconditioners for fluid-structure systems with large added-mass effect},
  Comput. Methods Appl. Mech. Engrg., 197 (2008), pp.~4216--4232.

\bibitem{BS08}
{\sc S.~Badia, A.~Quaini, and A.~Quarteroni}, {\em Splitting methods based on
  algebraic factorization for fluid-structure interaction}, SIAM J. Sci.
  Comput., 30 (2008), pp.~1778--1805.

\bibitem{ABXCC10}
{\sc A.~Barker and X.~Cai}, {\em Scalable parallel methods for monolithic
  coupling in fluid-structure interaction with application to blood flow
  modeling}, J. Comput. Phys.,  (2010), pp.~642--659.

\bibitem{JB08:00}
{\sc J.~Bonet and R.~Wood}, {\em Nonlinear Continuum Mechanics for Finite
  Element Analysis}, Cambridge University Press, New York, 2008.

\bibitem{Causin04}
{\sc P.~Causin, J.~Gerbeau, and F.~Nobile}, {\em Added-mass effect in the
  design of partitioned algorithms for fluid-structure problems}, Comput.
  Methods Appl. Mech. Engrg., 194 (2005), pp.~4506--4527.

\bibitem{EPFL11}
{\sc P.~Crosetto, S.~Deparis, G.~Fourestey, and A.~Quarteroni}, {\em Parallel
  algorithms for fluid-structure interaction problems in haemodynamics}, SIAM
  J. Sci. Comput., 33 (2011), pp.~1598--1622.

\bibitem{JD04:00}
{\sc J.~Donea, A.~Huerta, J.~Ponthot, and A.~Ferran}, {\em {Arbitrary
  Lagrangian-Eulerian methods}}, in The Encyclopedia of Computational
  Mechanics, E.~Stein, R.~Borst, and T.~Hughes, eds., vol.~1, Wiley\& Sons,
  Ltd, 2004, pp.~413--437.

\bibitem{EH99}
{\sc H.~Elman}, {\em Preconditioning for the steady-state {N}avier--{S}tokes
  equations with low viscosity}, SIAM J. Sci. Comput., 20 (1999),
  pp.~1299--1316.

\bibitem{LF99:00}
{\sc L.~Formaggia and F.~Nobile}, {\em A stability analysis for the arbitrary
  {L}agrangian {E}ulerian formulation with finite elements}, East-West J.
  Numer. Math., 7 (1999), pp.~105--132.

\bibitem{NME:NME3001}
{\sc M.~W. Gee, U.~Küttler, and W.~A. Wall}, {\em Truly monolithic algebraic
  multigrid for fluid–structure interaction}, Int. J. Numer. Meth. Engng., 85
  (2011), pp.~987--1016.

\bibitem{UH:02}
{\sc G.~Haase and U.~Langer}, {\em Modern Methods in Scientific Computing and
  Applications}, vol.~75 of NATO Science Series II. Mathematics, Physics and
  Chemistry, Kluwer Academic Press, Dordrecht, 2002, ch.~Multigrid Methods:
  From Geometrical to Algebraic Versions, pp.~103--154.

\bibitem{WH03}
{\sc W.~Hackbusch}, {\em Multi-Grid Methods and Applications}, Springer,
  Berlin, 2003.

\bibitem{Heil20041}
{\sc M.~Heil}, {\em An efficient solver for the fully coupled solution of
  large-displacement fluid-structure interaction problems}, Comput. Methods
  Appl. Mech. Engrg., 193 (2004), pp.~1--23.

\bibitem{GAH00:00}
{\sc G.~Holzapfel}, {\em Nonlinear Solid Mechanics: A Continuum Approach for
  Engineering}, John Wiley \& Sons, Chichester, 2000.

\bibitem{STurek06}
{\sc J.~Hron and S.~Turek}, {\em A monolithic {FEM/Multigrid} solver for an
  {ALE} formulation of fluid-structure interaction with applications in
  biomechanics}, in Fluid-Structure Interaction, H.-J. Bungartz and
  M.~Sch\"{a}fer, eds., vol.~53 of Lecture Notes in Computational Science and
  Engineering, Springer, 2006, pp.~146--170.

\bibitem{TH:81}
{\sc T.~Hughes, W.~Liu, and T.~Zimmermann}, {\em {Lagrangian-Eulerian finite
  element formulation for incompressible viscous flows}}, Comput. Methods Appl.
  Mech. Engrg., 29 (1981), pp.~329--349.

\bibitem{UJ:91}
{\sc M.~Jung and U.~Langer}, {\em Application of multilevel methods to
  practical problems}, Surv. Math. Ind., 1 (1991), pp.~217--257.

\bibitem{FK98:00}
{\sc F.~Kickinger}, {\em Algebraic multigrid for discrete elliptic second-order
  problems}, in Multigrid Methods V. Proceedings of the 5th European Multigrid
  conference (ed.~by W. Hackbush), Lecture Notes in Computational Sciences and
  Engineering, vol.~3, Springer, 1998, pp.~157--172.

\bibitem{KJMS13}
{\sc J.~Kraus and S.~Margenov}, {\em Robust Algebraic Multilevel Methods and
  Algorithms}, vol.~5 of Radon Series on Computational and Applied Mathematics,
  Walter de Gruyter, Berlin, New York, 2009.

\bibitem{UY14}
{\sc U.~Langer and H.~Yang}, {\em Numerical simulation of fluid-structure
  interaction problems with hyperelastic models: A monolithic approach},
  arXiv:1408.3737,  (2014).

\bibitem{CVULHY1300}
{\sc U.~Langer and H.~Yang}, {\em Partitioned solution algorithms for
  fluid-structure interaction problems with hyperelastic models}, J. Comput.
  Appl. Math., 276 (2015), pp.~47--61.

\bibitem{HJPMRO95}
{\sc H.-P. Moran and R.~Ohayon}, {\em Fluid-Structure Interaction: Applied
  Numerical Methods}, John Wiley \& Sons, 1995.

\bibitem{Muddle20127315}
{\sc R.~L. Muddle, M.~Mihajlovi\'{c}, and M.~Heil}, {\em An efficient
  preconditioner for monolithically-coupled large-displacement fluid-structure
  interaction problems with pseudo-solid mesh updates}, J. Comput. Phys., 231
  (2012), pp.~7315--7334.

\bibitem{NLA:NLA542}
{\sc Y.~Notay and P.~S. Vassilevski}, {\em {Recursive Krylov-based multigrid
  cycles}}, Numer. Linear Algebra Appl., 15 (2008), pp.~473--487.

\bibitem{QuarteroniVali:1999}
{\sc A.~Quarteroni and A.~Valli}, {\em Domain Decomposition Methods for Partial
  Differential Equations}, Oxfort Sciences Publications, 1999.

\bibitem{Razzaq20121156}
{\sc M.~Razzaq, H.~Damanik, J.~Hron, A.~Ouazzi, and S.~Turek}, {\em {FEM}
  multigrid techniques for fluid–structure interaction with application to
  hemodynamics}, Appl. Numer. Math., 62 (2012), pp.~1156--1170.

\bibitem{JWRKS87}
{\sc J.~W. Ruge and K.~St{\"u}ben}, {\em Algebraische mehrgittermethoden
  ({AMG})}, in Multigrid Methods, vol.~3 of Frontiers in Applied Mathematics,
  SIAM, Philadelphia, 1987, pp.~73--130.

\bibitem{Saad1993}
{\sc Y.~Saad}, {\em A flexible inner-outer preconditioned {GMRES} algorithm},
  SIAM J. Sci. Comput., 14 (1993), pp.~461--469.

\bibitem{YS03:00}
{\sc Y.~Saad}, {\em Iterative Methods for Sparse Linear Systems}, SIAM,
  Philadelphia, 2003.

\bibitem{Saad86}
{\sc Y.~Saad and M.~H. Schultz}, {\em {GMRES}: {A} generalized minimal residual
  algorithm for solving nonsymmetric linear systems}, SIAM J. Sci. Stat.
  Comput., 7 (1986), pp.~856--869.

\bibitem{JS97:00}
{\sc J.~Sch\"{o}berl}, {\em {NETGEN - An advancing front 2D/3D-mesh generator
  based on abstract rules}}, Comput Visual Sci, 1 (1997), pp.~41--52.

\bibitem{ATOW05}
{\sc A.~Toselli and O.~Widlund}, {\em Domain Decomposition Methods-Algorithms
  and Theory}, Springer, Heidelberg, 2005.

\bibitem{ST99:00}
{\sc S.~Turek}, {\em Efficient Solvers for Incompressible Flow Problems},
  Springer, Berlin, 1999.

\bibitem{PSV08:00}
{\sc P.~S. Vassilevski}, {\em Multilevel Block Factorization Preconditioners},
  Springer, Heidelberg, 2008.

\bibitem{WM04:00}
{\sc M.~Wabro}, {\em Coupled algebraic multigrid methods for the {O}seen
  problem}, Comput. Visual. Sci., 7 (2004), pp.~141--151.

\bibitem{WM06:00}
\leavevmode\vrule height 2pt depth -1.6pt width 23pt, {\em {AMGe}---coarsening
  strategies and application to the {O}seen equations}, SIAM J. Sci. Comput.,
  27 (2006), pp.~2077--2097.

\bibitem{Wick11:00}
{\sc T.~Wick}, {\em Fluid-structure interactions using different mesh motion
  techniques}, Comput. Structures, 89 (2011), pp.~1456--1467.

\bibitem{NLA:NLA1889}
{\sc T.~A. Wiesner, R.~S. Tuminaro, W.~A. Wall, and M.~W. Gee}, {\em {Multigrid
  transfers for nonsymmetric systems based on Schur complements and Galerkin
  projections}}, Numer. Linear Algebra Appl., 21 (2014), pp.~415--438.

\bibitem{YH11:00}
{\sc H.~Yang}, {\em Partitioned solvers for the fluid-structure interaction
  problems with a nearly incompressible elasticity model}, Comput. Visual.
  Sci., 14 (2011), pp.~227--247.

\bibitem{Yang11:00}
{\sc H.~Yang and W.~Zulehner}, {\em Numerical simulation of fluid-structure
  interaction problems on hybrid meshes with algebraic multigrid methods}, J.
  Comput. Appl. Math., 235 (2011), pp.~5367--5379.

\end{thebibliography}
\bibliographystyle{siam}

\end{document}